\documentclass[doublespacing]{elsart3-1}

\usepackage{amssymb}

\usepackage[english,francais]{babel}

\newtheorem{theorem}{Theorem}[section]

\newtheorem{e-proposition}[theorem]{Proposition}
\newtheorem{corollary}[theorem]{Corollary}
\newtheorem{e-definition}[theorem]{Definition\rm}

\renewcommand{\theequation}{\arabic{equation}}
\def\og{\leavevmode\raise.3ex\hbox{$\scriptscriptstyle\langle\!\langle$~}}
\def\fg{\leavevmode\raise.3ex\hbox{~$\!\scriptscriptstyle\,\rangle\!\rangle$}}

\journal{Comptes Rendus Mathematique}
\begin{document}
% place in the next line the header (rubrique) chosen for your article,
% if you know it (you can also have 2, format : Header1/Header2
\centerline{}
\begin{frontmatter}


% \author{Name\thanksref{label2}}
\selectlanguage{english}
\title{Unique continuation estimates for the Kolmogorov equation in the whole space}

% The [label1] can be suppressed if there is only one address for all authors

\selectlanguage{english}
\author[authorlabel1]{Yubiao Zhang},
\ead{yubiao\b{ }zhang@whu.edu.cn}
%\author[authorlabel2]{Author Name2}
%\ead{author.name2@email.address2}

\address[authorlabel1]{School of
Mathematics and Statistics, Wuhan University, Wuhan, 430072, China }
%\address[authorlabel2]{Address2}

% If you know the dates of reception, and acceptation you can put them now;
%  idem the name of the person presenting the Note

\medskip
%\begin{center}
%{\small Received *****; accepted after revision +++++\\
%Presented by £££££}
%\end{center}

\begin{abstract}
\selectlanguage{english}
% Text of abstract in English
We prove in this Note an observation estimate at one point in time for the
Kolmogorov equation in the whole space. Such estimate implies the observability and the null controllability
for the Kolmogorov equation with a control region which is sufficiently spread out throughout
the whole space.
%{\it To cite this article: A.
%Name1, A. Name2, C. R. Acad. Sci. Paris, Ser. I 340 (2005).}

\vskip 0.5\baselineskip

\selectlanguage{francais}
% Text of abstract in French
\noindent{\bf R\'esum\'e} \vskip 0.5\baselineskip \noindent
{\bf In\'egalit\'es de continuation unique pour l'\'equation de Kolmogorov dans l'espace tout entier. }
Nous montrons dans cette Note des in\'egalit\'es d'observation traduisant la continuation unique
pour l'\'equation de Kolmogorov d\'efinie sur l'espace tout entier.
%{\it Pour citer cet article~: A. Name1, A. Name2, C. R. Acad. Sci.
%Paris, Ser. I 340 (2005).}

\end{abstract}
\end{frontmatter}

% now the Version française abrégée, if it exists
\selectlanguage{francais}
\section*{Version fran\c{c}aise abr\'eg\'ee}
% Text of your Version française abrégée here.
% Note you do not need to repeat here equations that you use in the
% main text - for example 'voir (3)' is quite acceptable.

\selectlanguage{english}
% main text
\section{Introduction and the main result}
\label{}

\noindent Consider the following Kolmogorov equation in the whole space ($d\in\mathbb N^+$)
\begin{eqnarray}\label{K-eq}
\left\{\begin{array}{ll}
        (\partial_t+v\cdot\nabla_x-\Delta_v)g(t,x,v)=0, &(t,x,v)\in\mathbb R^+\times\mathbb R^d\times\mathbb R^d,\\
        g(0,x,v)=g_0(x,v),  &(x,v)\in\mathbb R^d\times\mathbb R^d.
       \end{array}
\right.
\end{eqnarray}
The well-posedness of the solution of (\ref{K-eq}) was proved in Propositions 2.1 and 2.2 in \cite{RM}. In \cite{RM}, the authors considered the following definition.\\

\begin{e-definition}[See Definition 1.1 in \cite{RM}]
An open set $O$ of $\mathbb R^n$ ($n\in\mathbb N^+$) is said to be an observability open set on the whole space $\mathbb R^n$ if there exist $\delta>0$ and $r>0$ such that
\begin{eqnarray}\label{geometry}
 \forall\,y\in\mathbb R^n,~\exists\,y^\prime\in O \mbox{ such that }
 B_{\mathbb R^n}(y^\prime,r)\subset O \mbox{ and } |y-y^\prime|\leq \delta.
\end{eqnarray}
Here $B_{\mathbb R^n}(x,r)$ denotes  an open ball in $\mathbb R^n$ of radius $r$ centered at $x$.\\

\end{e-definition}

From this definition, the authors in \cite{RM} proved the following estimate: Assume that  $\omega_x\subset\mathbb R^d$ and $\omega_v\subset\mathbb R^d$ both verifies the property (\ref{geometry}) with $n=d$. Then for all $T>0$, there exists $C>0$ so that for each $g_0\in L^2(\mathbb R^{2d})$, the solution of (\ref{K-eq}) satisfies that
\begin{eqnarray}\label{RM-observe}
 \|g|_{t=T}\|_{L^2(\mathbb R^{2d})}\leq C\|g\|_{L^2((0,T)\times\omega_x\times\omega_v)}.
\end{eqnarray}
In \cite{RM}, the proof of (\ref{RM-observe}) is based on a spectral inequality, a Carleman inequality with respect to the variable $v$ and a decay inequality for the Fourier transform of the solution of (\ref{K-eq}) with respect to the variable $x$. The geometric condition (\ref{geometry}) plays an important role in proving (\ref{RM-observe}). The authors in \cite{RM} pointed out the following fact: There exists an open set $\mathcal O$ of $\mathbb R^{2d}$, which is an observability open set in the whole $\mathbb R^{2d}$, and does not contain any cartesian product $\mathcal O_1\times \mathcal O_2$, where each $\mathcal O_1$ and $\mathcal O_2$ are both observability open sets in the whole space $\mathbb R^d$.

In this Note, when assume that $\omega\subset\mathbb R^{2d}$ verifies (\ref{geometry}) with $n=2d$, we  get a unique continuation estimate for the Kolmogorov equation. Such kind of estimate  has been studied in \cite{AEWZ} and \cite{PW}. Our proof combines  the spectral inequality given in \cite{RM} and  a decay inequality on the Fourier transform of the solution of (\ref{K-eq}) with respect to the variables $x$ and $v$. The main result is as follows.\\

%\vskip5pt
\begin{theorem}\label{Theorem-interpo-ineq}
Let $\omega\subset\mathbb R^{2d}$ be an observability open set on the whole space $\mathbb R^{2d}$. Then there exists  $C=C(\omega,d)>0$  so that for all $T>0$, $\alpha\in (0,1)$ and  $g_0\in L^2(\mathbb R^{2d})$, the solution of (\ref{K-eq}) satisfies that
\begin{eqnarray}\label{interpo-ineq}
 \big\|g|_{t=T}\big\|_{L^2(\mathbb R^{2d})}  \leq e^{\frac{C}{\alpha}(1+\frac{1}{T^3})}  \big\| g|_{t=T}\big\|_{L^2(\omega)}^{1-\alpha} \|g_0\|_{L^2(\mathbb R^{2d})}^{\alpha}.
\end{eqnarray}
\end{theorem}

%\vskip 10pt
%\begin{remark}
%The orders  of $\alpha$ and $T$ of the constant $e^{C(1+\frac{1}{\alpha T^3})}$ in (\ref{interpo-ineq}) are optimal.
%
%\end{remark}

\vskip 5pt
By a telescoping series method (see \cite[Theorem1.1]{PW}), a direct consequence  of  (\ref{interpo-ineq})  is the following observability estimate.\\
%\vskip5pt

\begin{corollary}\label{Corollary-measurable}
 Let $\omega\subset\mathbb R^{2d}$ be an observability open set on the whole space $\mathbb R^{2d}$. Let $T>0$ and $E\subset(0,T)$ be a measurable set of positive measure. Then there exists  $C_{obs}=C(\omega,d,T,E)>0$  so that  for each  $g_0\in L^2(\mathbb R^{2d})$, the solution of (\ref{K-eq}) verifies that
\begin{eqnarray}\label{interpo-ineq-meas}
 \big\|g|_{t=T}\big\|_{L^2(\mathbb R^{2d})}  \leq C_{obs}
 \int_E \|g(t,\cdot,\cdot)\|_{L^2(\omega)} \mathrm dt.
\end{eqnarray}
When $E=(0,T)$, $C_{obs}= e^{C(1+\frac{1}{T^3})}$ where $C$ only depends on $\omega$ and $d$.\\
\end{corollary}

Such observability estimate implies  by duality the null controllability  for the Kolmogorov equation.

\section{A spectral inequality}

%Let $\omega_x\subset\mathbb R^{d}$ be an observability open set on the whole space $\mathbb R^d$. Then there exists a constant $C>0$ such that
%\begin{eqnarray*}
% \|f\|_{L^2(\mathbb R^d)} \leq e^{C(N+1)} \|f\|_{L^2(\omega_x)}
%\end{eqnarray*}
%for $N\geq0$ and $f\in L^2(\mathbb R^d)$ such that supp$(\hat f)\subset\overline{B_{\mathbb R^d}(0,N)}$, the closed ball of radius $N$ and center $0$.

\noindent
%In this and the following sections, when $f\in L^2(\mathbb R^{2d})$, we denote by $\hat f$  the Fourier transform of $f$ with respect to the variables $x$ and $v$.
 The following spectral inequality plays a key role to deduce the estimate (\ref{interpo-ineq}). Here $\hat f$ denotes the Fourier transform  of $f$.\\

 \begin{theorem}[See Theorem 1.2 in \cite{RM}]\label{theorem-spectral}
 Let $\omega\subset\mathbb R^{2d}$ be an observability open set on the whole space $\mathbb R^{2d}$. Then there exists  $C=C(\omega,d)>0$ such that  for all $N> 0$, every $f\in L^2(\mathbb R^{2d})$ verifies that
\begin{eqnarray}\label{spec-ineq-wholespace}
 \int_{|\xi|\leq N} |\hat f(\xi)|^2 \mathrm d\xi \leq e^{C(1+N)} \int_{\omega} \Big|\int_{|\xi|\leq N}  \hat f(\xi) e^{i x\xi} \mathrm d\xi \Big|^2 \mathrm dx.
\end{eqnarray}
 \end{theorem}

\vskip 5pt

 We mention that, for smooth  compact and connected Riemannian manifold $M$ with metric $g$ and boundary $\partial M$, the following inequality was obtained in \cite{LZ}: Let $\omega\subset M$ be an open nonempty subset. There exists $C>0$ such that the Laplace-Beltrami operator $-\Delta_g$ on $M$ satisfies that
\begin{eqnarray}\label{spectral-M}
 \|u\|_{L^2(M)} \leq C e^{C\sqrt{\lambda}} \| u\|_{L^2(\omega)}  \mbox{ for all } u\in\mbox{span} \{e_j;\,\lambda_j\leq \lambda\},
\end{eqnarray}
where  $\{\lambda_j\}$ and $\{e_j\}$ are the eigenvalues and the corresponding eigenvectors of $-\Delta_g$ with the zero Dirichlet boundary condition.  Based on this type of inequality (\ref{spectral-M}), a similar estimate to (\ref{interpo-ineq}) was obtained for the heat equation in a bounded domain (see \cite[Theorem 6]{AEWZ}). The strategy in this Note also works for the heat equation in the whole space. This can be compared with \cite{Miller}, where $M$ is non-compact with  a Ricci curvature  bounded below. The author in \cite{Miller} proves that, under an interpolation inequality in \cite[(6) on Page 40]{Miller}, (\ref{geometry}) implies the spectral inequality (\ref{spec-ineq-wholespace}), which yields the observability for the heat equation in $M$.

%For non-compact manifold $M$, we mention the work in \cite{Miller}, some sufficient and some necessary geometric conditions of (\ref{spectral-M}) are discussed.

%%\begin{proof}
%Let $N>0$ and $f\in L^2(\mathbb R^{2d})$. Consider the function
%\begin{eqnarray*}
% g(\xi)\triangleq \chi_{\{|\xi|\leq N\}} \hat f(\xi),~\forall\,\xi\in\mathbb R^{2d}.
%\end{eqnarray*}
%Then one can directly check that
%\begin{eqnarray*}
%  \int_{|\xi|\leq N} |\hat f(\xi)|^2 \mathrm d\xi =\int_{\mathbb R^{2d}_\xi} |g(\xi)|^2 \mathrm d\xi = \frac{1}{(2\pi)^{2d}} \int_{\mathbb R^{2d}_x} |\mathcal F^{-1}(g)(x)|^2 \mathrm dx
%\end{eqnarray*}
%and
%\begin{eqnarray*}
% \int_{\omega} \Big|\int_{|\xi|\leq N}  \hat f(\xi) e^{i x\xi} \mathrm d\xi \Big|^2 \mathrm dx=\int_{\omega} |\mathcal F^{-1}(g)(x)|^2 \mathrm dx,
%\end{eqnarray*}
%where $\mathcal F^{-1}$ denotes the inverse Fourier transform. The above, as well as the spectral inequality proved by J. Rousseau and I. Moyano, yields (\ref{spec-ineq-wholespace}).
%
%%\end{proof}

\section{A decay inequality}

\noindent We apply the Fourier transform, with respect to the variables $x$ and $v$, to Equation (\ref{K-eq}). Then we get the following equation in the corresponding frequency space
\begin{eqnarray}\label{K-eq-adjoint}
\left\{\begin{array}{ll}
        (\partial_t-\xi\cdot\nabla_\eta+|\eta|^2)\hat g(t,\xi,\eta)=0, &(t,\xi,\eta)\in\mathbb R^+\times\mathbb R^d\times\mathbb R^d,\\
        \hat g(0,\xi,\eta)=\hat g_0(\xi,\eta),  &(\xi,\eta)\in\mathbb R^d\times\mathbb R^d.
       \end{array}
\right.
\end{eqnarray}
The solution of (\ref{K-eq-adjoint}) has an explicit representation, which has been obtained in \cite[Section 7.6, Pages 210-211]{H}. Based on this, we get a decay estimate for the Kolmogorov equation as follows.\\

\begin{e-proposition}\label{lemma-decay}
There exist $C>0$ and $C^\prime=C^\prime(d)>0$ such that for all $N$, $T>0$ and each $g_0\in L^2(\mathbb R^{2d})$, the solution of (\ref{K-eq-adjoint}) verifies that
\begin{eqnarray}\label{A-decay-condition}
 \int_{|(\xi,\eta)|> N} |\hat g(T,\xi,\eta)|^2 \mathrm d\xi\mathrm d\eta \leq  e^{C^\prime-CN^2 \min\{T,T^3\} }  \int_{\mathbb R_x^{d}\times\mathbb R_v^d} |g_0(x,v)|^2 \mathrm d x\mathrm d v.
\end{eqnarray}
\end{e-proposition}

\vskip5pt
\noindent\textbf{Proof.} Let $ g$ be a solution of  (\ref{K-eq-adjoint}). One can directly compute   that
\begin{eqnarray*}
 \hat g(t,\xi,\eta)= \hat g_0(\xi,\eta+\xi t)  \exp\big( -|\eta|^2 t-\eta\cdot\xi t^2-|\xi|^2 t^3/3 \big),~\forall\,(t,\xi,\eta)\in\mathbb R^+\times\mathbb R^d\times\mathbb R^d.
\end{eqnarray*}
This yields that for all $(t,\xi,\eta)\in\mathbb R^+\times\mathbb R^d\times\mathbb R^d$,
\begin{eqnarray*}
 |\hat g(t,\xi,\eta)|
 &\leq& |\hat g_0(\xi,\eta+\xi t)|  \exp\big[ -(|\eta|^2 +|\xi|^2 )\min\{t, t^3\}/30 \big].
\end{eqnarray*}
From this, we see that for all $N$, $T>0$,
\begin{eqnarray*}
\int_{|(\xi,\eta)|> N} |\hat g(T,\xi,\eta)|^2 \mathrm d\xi\mathrm d\eta&\leq& \exp\big( -N^2\min\{T, T^3\}/15 \big)  \int_{\mathbb R_\xi^{d}\times\mathbb R_\eta^d} |\hat g_0(\xi,\eta)|^2 \mathrm d\xi\mathrm d\eta,
\end{eqnarray*}
which leads to (\ref{A-decay-condition}). This ends the proof.
\begin{flushright}
  \qed
\end{flushright}

\section{Proofs of Theorem \ref{Theorem-interpo-ineq} and Corollary \ref{Corollary-measurable}}

\noindent In this section, we first  prove Theorem \ref{Theorem-interpo-ineq} by combining Theorem \ref{theorem-spectral} and Proposition \ref{lemma-decay} as follows.\\

 \noindent\textbf{Proof of Theorem \ref{Theorem-interpo-ineq}.} Let   $g$ be the solution of Equation (\ref{K-eq}) with the initial data $g_0\in L^2(\mathbb R^{2d})$. For each $N>0$, write
\begin{eqnarray*}
\hat g(t,\xi,\eta)&=&\chi_{B_N}(\xi,\eta) \hat g(t,\xi,\eta) +
\chi_{B_N^c}(\xi,\eta) \hat g(t,\xi,\eta)
,~\forall\,(t,\xi,\eta)\in\mathbb R^+\times\mathbb R^d\times\mathbb R^d,
\end{eqnarray*}
where $\chi_{B_N}$ and  $\chi_{B_N^c}$ denote the characteristic functions of the set $B_N\triangleq\big\{(\xi,\eta)\in\mathbb R^{2d};\,|(\xi,\eta)|\leq N \big\}$ and its complement, respectively. Let $T>0$. We observe that for all $N>0$,
\begin{eqnarray}\label{inetpo-1}
 (2\pi)^{d}\|g|_{t=T}\|_{L^2(\mathbb R^{2d})} = \|\hat g|_{t=T}\|_{L^2(\mathbb R^{2d})}\leq \|\chi_{B_N} \hat g|_{t=T}\|_{L^2(\mathbb R^{2d})}+\|\chi_{B_N^c} \hat g|_{t=T}\|_{L^2(\mathbb R^{2d})}.~~~~
\end{eqnarray}
On one hand, we  apply (\ref{spec-ineq-wholespace}) to $g$ to get the existence of a positive constant $C_1=C_1(\omega,d)$ so that for all $N>0$,
\begin{eqnarray}\label{inetpo-1-1}
  \int_{B_N} |\hat g(T,\xi,\eta)|^2 \mathrm d\xi\mathrm d\eta
  %\nonumber\\
%   &\leq& e^{2C_1(N+1)} \int_\omega \Big\|\int_{B_N}\hat g(T,\xi,\eta)e^{i(x\xi+v\eta)} \mathrm d\xi\mathrm d\eta \Big\|^2 \mathrm d x\mathrm d v \nonumber\\
 &\leq& e^{2C_1(N+1)}  \Big[ \int_\omega \big|\int_{\mathbb R^d_\xi \times \mathbb R^d_\eta}  \hat g(T,\xi,\eta)e^{i(x\cdot\xi+v\cdot\eta)} \mathrm d\xi\mathrm d\eta \big|^2 \mathrm d x\mathrm d v    \nonumber\\
 & & +  \int_{\mathbb R^d_x \times \mathbb R^d_v} \big|\int_{B_N^c}\hat g(T,\xi,\eta)e^{i(x\cdot\xi+v\cdot\eta)} \mathrm d\xi\mathrm d\eta \big|^2 \mathrm d x\mathrm d v   \Big].
\end{eqnarray}
On the other hand, let $f(\xi,\eta)\triangleq\chi_{B_N^c}(\xi,\eta)\hat g(T,\xi,\eta)$, $(\xi,\eta)\in  \mathbb R^d_\xi \times \mathbb R^d_\eta$. It follows from the inverse Fourier transform formula that $\int f(\xi,\eta)e^{i(x\cdot\xi+v\cdot\eta)} \mathrm d\xi\mathrm d\eta$ is the inverse Fourier transform of $f$. Then
\begin{eqnarray}\label{inetpo-1-2}
 \frac{1}{(2\pi)^{2d}} \int_{\mathbb R^d_x \times \mathbb R^d_v} \big|\int_{B_N^c}\hat g(T,\xi,\eta)e^{i(x\cdot\xi+v\cdot\eta)} \mathrm d\xi\mathrm d\eta \big|^2 \mathrm d x\mathrm d v
 &=& \frac{1}{(2\pi)^{2d}} \int_{\mathbb R^d_x \times \mathbb R^d_v} \big|\int_{\mathbb R^d_\xi \times \mathbb R^d_\eta} f(\xi,\eta)e^{i(x\cdot\xi+v\cdot\eta)} \mathrm d\xi\mathrm d\eta \big|^2 \mathrm d x\mathrm d v  \nonumber\\
 &=& \int_{\mathbb R_\xi^{d}\times\mathbb R_\eta^d} |f(\xi,\eta)|^2 \mathrm d\xi\mathrm d\eta=\int_{B_N^c} |\hat g(T,\xi,\eta)|^2 \mathrm d\xi\mathrm d\eta.
\end{eqnarray}
Meanwhile, we apply (\ref{A-decay-condition}) to $g$ to obtain that there exist $C_2>0$ and $C_3=C_3(d)>0$ so that for all $N>0$,
\begin{eqnarray}\label{inetpo-1-3}
 \int_{B_N^c} |\hat g(T,\xi,\eta)|^2 \mathrm d\xi\mathrm d\eta \leq e^{2[C_3-C_2N^2\min\{T,T^3\}] }  \int_{\mathbb R_x^{d}\times\mathbb R_v^d} |g_0(x,v)|^2 \mathrm d x\mathrm d v.
\end{eqnarray}
Write $T_3^1\triangleq\min\{T,T^3\}$. By the inverse Fourier transform formula, we see from (\ref{inetpo-1})-(\ref{inetpo-1-3}) that  for all $N>0$,
\begin{eqnarray}\label{inetpo-3}
 \|g|_{t=T}\|_{L^2(\mathbb R^{2d})} &\leq& e^{C_1(N+1)}\|g|_{t=T}\|_{L^2(\omega)}  + 2e^{C_1(N+1)+C_3-C_2N^2 T_3^1} \|g_0\|_{L^2(\mathbb R^{2d})}.~~~
\end{eqnarray}
Let $\alpha\in (0,1)$. We set $ k(\alpha)\triangleq \alpha/(1-\alpha)$.
Then we have that for all $N>0$,
\begin{eqnarray*}
 C_1N &\leq &  \frac{C_1^2}{2k(\alpha)C_2 T_3^1}  + k(\alpha)\frac{C_2 N^2 T_3^1}{2}
\mbox{ and }
  C_1N-C_2N^2 T_3^1  \leq \frac{C_1^2}{2C_2 T_3^1}-\frac{C_2 N^2 T_3^1}{2}.
\end{eqnarray*}
These, together with (\ref{inetpo-3}), yield that for all $\varepsilon\in(0,1)$,
\begin{eqnarray}\label{inetpo-5}
 \|g|_{t=T}\|_{L^2(\mathbb R^{2d})}
 \leq \widetilde C_1 \Big[ \varepsilon^{-k(\alpha)}  \| g|_{t=T}\|_{L^2(\omega)}+
 \varepsilon \|g_0\|_{L^2(\mathbb R^{2d})} \Big],~~~
\end{eqnarray}
where
$$
 \widetilde C_1\triangleq \max\Big\{ e^{C_1+\frac{C_1^2}{2k(\alpha) C_2T_3^1}},  2e^{C_1+C_3+\frac{C_1^2}{2C_2T_3^1}} \Big\}\leq 2e^{\frac{(C_1+C_2+C_3)^2}{C_2\alpha}(1+\frac{1}{T^3})}.
 $$
Since $\|g|_{t=T}\|_{L^2(\mathbb R^{2d})}\leq \|g_0\|_{L^2(\mathbb R^{2d})}$, the  minimization of the right side of (\ref{inetpo-5}), with respect to the variable $\varepsilon$ over $\mathbb R^+$, leads to (\ref{interpo-ineq}). This completes the proof.
\begin{flushright}
  \qed
\end{flushright}

\vskip5pt
\noindent We next use  the telescoping series method to deduce the Corollary \ref{Corollary-measurable} from Theorem \ref{Theorem-interpo-ineq}.\\

\noindent\textbf{Proof of Corollary \ref{Corollary-measurable}.}  Let   $g$ be the solution of Equation (\ref{K-eq}) with the initial data $g_0\in L^2(\mathbb R^{2d})$.  We take $\alpha=1/2$ in (\ref{interpo-ineq}) and then see from  the Young inequality that there exists $C_1=C_1(\omega,d)>0$ so that
\begin{eqnarray*}
 \|g|_{t=T}\|     \leq    \frac{1}{\varepsilon} e^{C_1(1+\frac{1}{T^3})}  \|g(T,\cdot,\cdot)\|_{L^2(\omega)}
 +\varepsilon \|g|_{t=0}\|,~\forall\,\varepsilon>0.
\end{eqnarray*}
Generally, for each $0<t_1<t_2$, we have that
\begin{eqnarray}\label{meas-1}
 \|g|_{t=t_2}\|     \leq    \frac{1}{\varepsilon} e^{C_1[1+\frac{1}{(t_2-t_1)^3}]}  \|g(t_2,\cdot,\cdot)\|_{L^2(\omega)}
 +\varepsilon \|g|_{t=t_1}\|,~\forall\,\varepsilon>0.
\end{eqnarray}
Let $l$ be a Lebesgue density point of $E$. Then by \cite[Proposition 2.1]{PW}, we know that for each $\lambda\in (1/\sqrt[6]{2},1)$, there exists a sequence $\{l_m\}\subset(l,T)$ so that for each $m\in\mathbb N^+$,
\begin{eqnarray}\label{meas-2}
  l_m-l=\lambda^{m-1} (l_1-l)   \mbox{ and }
   3|E\cap(l_{m+1},l_m)|\geq  |l_{m+1}-l_m|.
\end{eqnarray}
Take a $m\in\mathbb N^+$ and let  $0<l_{m+2}<l_{m+1}\leq s<l_m<T$. Since $\|g|_{t=l_{m}}\|\leq \|g|_{t=s}\|$ and $l_{m+1}-l_{m+2} \leq s-l_{m+2}$, we apply (\ref{meas-1}), where $t_1=l_{m+2}$ and $t_2=s$, to get that
\begin{eqnarray*}
  \|g|_{t=l_{m}}\|\leq   \frac{1}{\varepsilon} e^{C_1[1+\frac{1}{(l_{m+1}-l_{m+2})^3}]}   \|g(s,\cdot,\cdot)\|_{L^2(\omega)}
 +\varepsilon \|g|_{t=l_{m+2}}\|,~\forall\,\varepsilon>0.
\end{eqnarray*}
By integrating both sides over $E\cap(l_{m+1},l_m)$ in the above inequality, we know that
%\begin{eqnarray*}
% |E\cap(l_{m+1},l_m)|\|g(l_m,\cdot,\cdot)\|  &\leq&   \frac{1}{\varepsilon} e^{C_1[1+\frac{1}{(l_m-l_{m+2})^3}]} \int_{E\cap(l_{m+1},l_m)}\|\chi_\omega g(t,\cdot,\cdot)\|\mathrm dt \nonumber\\
% & &  +\varepsilon|E\cap(l_{m+1},l_m)| \|g(l_{m+2},\cdot,\cdot)\|,~\forall\,\varepsilon>0.
%\end{eqnarray*}
%Therefore, for all $m\in\mathbb N^+$, we obtain that
\begin{eqnarray}\label{meas-3}
 & &  \Big(\varepsilon |E\cap(l_{m+1},l_m)| e^{-\frac{C_1}{(l_{m+1}-l_{m+2})^3}}\Big)   \|g|_{t=l_{m}}\| -\Big(\varepsilon^2|E\cap(l_{m+1},l_m)|e^{-\frac{C_1}{(l_{m+1}-l_{m+2})^3}}\Big)   \|g|_{t=l_{m+2}}\| \nonumber\\
 &\leq& e^{C_1} \int_{E\cap(l_{m+1},l_m)}\|g(s,\cdot,\cdot)\|_{L^2(\omega)}\mathrm ds,~\forall\,\varepsilon>0.
\end{eqnarray}
Meanwhile, we know from (\ref{meas-2}) that
\begin{eqnarray*}
 3|E\cap(l_{m+1},l_m)|\geq  |l_{m+1}-l_m|\geq  e^{- \frac{1}{|l_{m+1}-l_m|} }
  \geq e^{ -\frac{\lambda^3(l_1-l_2)^2}{(l_{m+1}-l_{m+2})^3} },~\forall\, m\in\mathbb N^+.
\end{eqnarray*}
Since $l_{m}-l_{m+2}=(1+\frac{1}{\lambda})(l_{m+1}-l_{m+2})$, the above, as well as (\ref{meas-3}), yields that for all $m\in\mathbb N^+$ and $\varepsilon>0$,
\begin{eqnarray}\label{meas-4}
   \varepsilon e^{-\frac{C_2}{(l_{m}-l_{m+2})^3}}\|g|_{t=l_{m}}\| -\varepsilon^2e^{-\frac{C_2}{(l_{m}-l_{m+2})^3}}\|g|_{t=l_{m+2}}\| \leq  3 e^{C_1} \int_{E\cap(l_{m+1},l_m)} \|g(s,\cdot,\cdot)\|_{L^2(\omega)}   \mathrm ds,
\end{eqnarray}
where $C_2=(1+\frac{1}{\lambda})^3[C_1+\lambda^3(l_1-l_2)^2]$. Let $\beta\triangleq\frac{\lambda^6}{2\lambda^6-1}$ $(>0)$ and $\varepsilon=e^{-\frac{(\beta-1)C_2}{(l_{m}-l_{m+2})^3}}$.
Since $\lambda^2(l_{m}-l_{m+2})=l_{m+2}-l_{m+4}$, $\forall\,m\in\mathbb N^+$, it follows from (\ref{meas-4}) that
\begin{eqnarray*}
  & &  e^{-\frac{\beta C_2}{(l_{m}-l_{m+2})^3}}\|g|_{t=l_{m}}\| -  e^{-\frac{\beta C_2 }{(l_{m+2}-l_{m+4})^3}}\|g|_{t=l_{m+2}}\|
  \leq 3e^{C_1} \int_{E\cap(l_{m+1},l_m)}  \|g(s,\cdot,\cdot)\|_{L^2(\omega)}   \mathrm ds.
\end{eqnarray*}
We deduce from this that
\begin{eqnarray*}
 e^{-\frac{\beta C_2}{(l_1-l_{3})^3}}\|g|_{t=l_1}\| &=& \sum_{k=0}^\infty  \Big[ e^{-\frac{\beta C_2}{(l_{2k+1}-l_{2k+3})^3}}\|g|_{t=l_{2k+1}}\| -  e^{-\frac{\beta C_2}{(l_{2k+3}-l_{2k+5})^3}}\|g|_{t=l_{2k+3}}\|  \Big] \nonumber\\
 &\leq& \sum_{k=0}^\infty 3e^{C_1} \int_{E\cap(l_{2k+3},l_{2k+1})}\|g(s,\cdot,\cdot)\|_{L^2(\omega)}\mathrm ds \leq 3e^{C_1}\int_{E\cap(l,l_1)}  \,\|g(s,\cdot,\cdot)\|_{L^2(\omega)}\,\mathrm ds.
\end{eqnarray*}
Since $\|g|_{t=T}\|\leq \|g|_{t=l_1}\|$, the above implies that
\begin{eqnarray*}
 \|g|_{t=T}\| \leq 3e^{C_1+\frac{\beta C_2 }{(l_1-l_{3})^3}}\int_{E}  \,\|g(s,\cdot,\cdot)\|_{L^2(\omega)}\,\mathrm ds.
\end{eqnarray*}
This proves (\ref{interpo-ineq-meas}). Especially, when $E=(0,T)$, we can take $l_1=T$ and $l_3=T/4$. We end the proof.
\begin{flushright}
  \qed
\end{flushright}

% etc, etc

% \appendix

% \section{}
% \label{}

% The Acknowledgements are an un-numbered section
%\section*{Acknowledgements}
\vskip 5pt
 \textbf{Acknowledgements.} The author gratefully thanks Professor Kim Dang Phung for  discussing and his valuable suggestions. Also, the author would like to thank Can Zhang for his help.

\end{document}